\documentclass[a4paper]{amsart}
\usepackage{graphicx}
\usepackage{amsmath}
\usepackage{amssymb}
\usepackage{oldgerm}
\usepackage[all]{xy}
\usepackage{amscd}
\usepackage{color}

\begin{document}

\newcommand{\blue}{\color{blue}}
\newcommand{\Hom}{\operatorname{Hom}\nolimits}
\newcommand{\uHom}{\operatorname{\underline{Hom}}\nolimits}
\newcommand{\End}{\operatorname{End}\nolimits}
\renewcommand{\Im}{\operatorname{Im}\nolimits}
\renewcommand{\mod}{\operatorname{mod}\nolimits}
\newcommand{\stmod}{\operatorname{\underline{mod}}\nolimits}
\newcommand{\Ker}{\operatorname{Ker}\nolimits}
\newcommand{\gr}{\operatorname{gr}\nolimits}
\newcommand{\Coker}{\operatorname{Coker}\nolimits}
\newcommand{\Coim}{{\operatorname{Coim}\nolimits}}
\newcommand{\Pd}{\operatorname{pd}\nolimits}
\newcommand{\gldim}{\operatorname{gldim}\nolimits}
\newcommand{\repdim}{\operatorname{repdim}\nolimits}
\newcommand{\Ht}{\operatorname{ht}\nolimits}
                                                              
\newcommand{\ra}{\operatorname{\mathfrak{r}}\nolimits}
\newcommand{\La}{\Lambda}
\newcommand{\Ga}{\Gamma}
\newcommand{\Z}{\operatorname{Z}\nolimits}
\newcommand{\cx}{\operatorname{cx}\nolimits}
\newcommand{\op}{\operatorname{op}\nolimits}
\newcommand{\V}{\operatorname{V}\nolimits}
\newcommand{\A}{\operatorname{A}\nolimits}
\newcommand{\T}{\operatorname{T}\nolimits}
\newcommand{\cidim}{\operatorname{CI-dim}\nolimits}
\newcommand{\e}{\operatorname{e}\nolimits}
\newcommand{\Lae}{\operatorname{\Lambda^{\e}}\nolimits}
\newcommand{\Ae}{\operatorname{A^{\e}}\nolimits}
\newcommand{\lcm}{\operatorname{lcm}\nolimits}
\newcommand{\cha}{\operatorname{char}\nolimits}
\newcommand{\add}{\operatorname{add}\nolimits}
\newcommand{\thick}{\operatorname{thick}\nolimits}

\newcommand{\Alt}{\mathcal{A}}
\newcommand{\cH}{\mathcal{H}}
\newcommand{\cS}{\mathcal{S}}
\newcommand{\hg}{\hat{G}}
\newcommand{\hG}{\hat{\Gamma}}

\newcommand{\vf}{\varphi}

\newcommand{\ul}{\underline}
\newcommand{\ve}{\varepsilon}
\newcommand{\ol}{\overline}

\newtheorem{theorem}{Theorem}[section]
\newtheorem{acknowledgement}[theorem]{Acknowledgement}
\newtheorem{ass}[theorem]{Assumption}
\newtheorem{expl}[theorem]{Example}
\newtheorem{case}[theorem]{Case}
\newtheorem{claim}[theorem]{Claim}
\newtheorem{conclusion}[theorem]{Conclusion}
\newtheorem{condition}[theorem]{Condition}
\newtheorem{conjecture}[theorem]{Conjecture}
\newtheorem{corollary}[theorem]{Corollary}
\newtheorem{criterion}[theorem]{Criterion}
\newtheorem{lemma}[theorem]{Lemma}
\newtheorem{prop}[theorem]{Proposition}
 \theoremstyle{definition}
\newtheorem*{definition}{Definition}
\theoremstyle{definition}
\newtheorem*{question}{Question}
\theoremstyle{definition}
\newtheorem*{questions}{Questions}
\theoremstyle{definition}
\newtheorem*{example}{Example}
\theoremstyle{definition}
\newtheorem*{examples}{Examples}
\newtheorem{exercise}[theorem]{Exercise}
\newtheorem{notation}[theorem]{Notation}
\theoremstyle{definition}
\newtheorem*{problem}{Problem}
                                                           
\theoremstyle{remark}
\newtheorem*{remark}{Remark}
\theoremstyle{remark}
\newtheorem*{remarks}{Remarks}
\theoremstyle{definition}
\newtheorem*{solution}{Solution}
\theoremstyle{definition}
\newtheorem*{assumption}{Assumption}
\newtheorem{summary}[theorem]{Summary}


\title{Cohomology of some local selfinjective algebras}
\author{Karin Erdmann}
\address{Karin Erdmann \\ Mathematical Institute \\ ROQ  \\ Oxford OX2 6GG \\ United Kingdom}
\email{erdmann@maths.ox.ac.uk}

\date{\today}



\begin{abstract}We show that the cohomology for some 2-generated selfinjective algebras is finitely 
generated. This applies in particular to the algebras $A5(\beta)$ for $\beta\neq 0$ in the classification 
of connected Hopf algebras of dimension $p^3$ over characteristic $p$, by Nguyen-Wang-Wang.
\end{abstract}

\maketitle

\section{Introduction}

Assume $A$ is a finite-dimensional Hopf algebra, this has a trivial module $k$, defined using the counit, and the cohomology of $A$ is the algebra
$H^*(A)= {\rm Ext}^*_{A}(k, k)$. 
More generally, a selfinjective algebra may have a trivial module, and then there is a cohomology. One would like to know
whether this cohomology is finitely generated; this is known for some
classes of algebras, but is open in general.

In this note we show that for some selfinjective local algebras
the cohomology is finitely generated. The novelty is that
we only use some identities for generators, rather than a complete
presentation of the algebra. Moreover, the proofs are completely
elementary.  This applies in particular to 
2-generated q-complete intersections as studied in \cite{BE}, \cite{EHF}, \cite{O}, and also 
to group algebras of 2-generated finite p-groups over fields of prime characteristic (though for 
group algebras, cohomology is known to be finitely generated \cite{E}.)

This is motivated by results in 
\cite{NWW} and \cite{ESW}. 
The work in \cite{NWW} classifies connected Hopf algebras of dimension $p^3$
over characteristic $p$, the result is a list of
24 families of algebras. These algebras and their
cohomology are analysed in \cite{ESW}. Perhaps surprisingly, the
algebra structures occuring are quite divers; by exploiting different
ideas, in all but one family 
it could be proved that the cohomology is finitely generated. 
The only case  left was the algebra with label $A5$ where char$(k)\geq 3$, with presentation
$$k\langle z, y\rangle /([y,z]^p, \ y^p, \ [[z,y]y],\  [[z, y], z],\
z^p + y[y,z]^{p-1} - \beta [y, z])
$$
with $\beta \in k$.   
For $\beta\neq 0$
we had mentioned in \cite{ESW} that its cohomology is finitely generated.
Our general result of this note shows that this  is indeed the case.
When $\beta=0$ the situation is different.

In Section 2 we describe the input. Section 3 computes the first few syzygies and discusses
further input.  In Section 4 we compute an explicit minimal projective resolution of $k$, and in 
Section 5 we show that the cohomology is finitely generated. The final section deals with the
case of the algebra $A5$ for $\beta \neq 0$.

Our proof only gives finite generation, but does not
compute a presentation of the cohomology. Such a presentation depends on the
input, as will be clear from illustrations in Sections 5 and 6.

\bigskip

\section{The input}\label{sec:2}

Let 
$\La$ 
be a local selfinjective algebra, of the form
$$\La  = k\langle x, y\rangle /I,$$ 
where $x$ and $y$ are independent generators of the radical $J=J(\La)$ of $\La$.
This has a trivial module, $k=\La/J$, and the cohomology of $\La$ is the algebra  
${\rm Ext}^*_{\La}(k, k)$. 

We impose conditions on
the ideal $I$, see below.
This includes several classes of algebras: group algebras of abelian $p$-groups of rank 2 over fields of characteristic $p$, also
quantum complete intersections.
As well this deals with 
the algebra $A5$ above when  $\beta\neq 0$ and char$(k)\geq 3$. 

\medskip

\begin{ass}\label{ass:2.1} \normalfont
We impose three  conditions on $I$, the first two are as follows.\\
(1) There is a minimal relation $\psi_1x + \psi_2y=0$ with $\psi_1, \psi_2$
independent modulo $J^2$. \\
(2)  There are minimal relations $\sigma_1x + \sigma_2y=0$ and $\theta_1x + \theta_2y=0$. Moreover
if $\sigma:= (\sigma_1, \sigma_2)$ and $\theta:= (\theta_1, \theta_2)$ then $\sigma$ and $\theta$ are
independent modulo $\La\psi$ where $\psi := (\psi_1, \psi_2)$. 
For example, such relations might be given via $x^n=0$ and $y^m=0$. 
\end{ass}

\bigskip 
The third condition will be formulated in Section \ref{sec:3}, in \ref{ass:3.4}.
\bigskip

\begin{expl} \label{ex:2.2} \normalfont Let 
$$A = \langle x, y\rangle /(x^n, y^m, xy-qyx)$$
 with  $0\neq q\in k$. This includes group algebras of
rank 2 finite abelian $p$-groups (when $k$ has characteristic $p$ and $q=1$, and $m, n$ are
powers of $p$) but also q-complete intersections when $q\neq 1$ as studied for example in \cite{BE}. 
 These algebras satisfy (1) and (2), taking
$$\psi = (-qy, x), \ \  \sigma = (x^{n-1}, 0), \ \ \theta = (0, y^{m-1}).
$$  
We use this algebra $A$ as an example throughout.  
\end{expl}

\bigskip

\section{The first few syzygies of $k$}\label{sec:3}

We work with left modules, and we write maps to the right.
First,  $\Omega(k) = \La x + \La y$.
Next, we take the second syzygy of $k$ in the form
$$\Omega^2(k) = \{ (a, b)\in \La^2 \mid ax + by = 0\}.
$$ 
By conditions  (1) and (2) of Assumption \ref{ass:2.1}, 
this module contains the independent generators:
$$\sigma = (\sigma_1, \sigma_2), \ \psi = (\psi_1, \psi_2), \ 
\theta = (\theta_1, \theta_2).
$$
The following will show that these in fact generate all of $\Omega^2(k)$.
Let $N:= \dim \La$.

\bigskip

\begin{lemma}\label{lem:3.1} \ We have  
	$\Omega^2(k) = \langle \sigma, \psi, \theta \rangle$ 
	as a left $\La$-module.
\end{lemma}

{\it Proof}  By conditions (1) and (2)  we have one inclusion, and we observe that $\dim \Omega^2(k) = N +1$. Now, since
the two components of $\psi$ are generators for $J$, and since $A$ is selfinjective, we have 
$\La\psi  \cong \Omega^{-1}(k)$ as a $\La$-module, in particular it has dimension
$N-1$. By assumption the elements $\sigma$ and $\theta$ are 
independent modulo $A\psi$. Hence $k\sigma + k\theta + \La\psi$ has dimension $n+1$, and  the claim follows.
$\Box$

\medskip

\begin{remark} \label{rem:3.1} For future reference, 
	the proof shows the following: The submodule
	$\langle \sigma, \psi \rangle$ of $\Omega^2(k)$ has dimension $N$, 
	and so does the submodule $\langle \psi, \theta \rangle$.
\end{remark}

\bigskip

We have the following immediate consequence.

\begin{lemma}\label{lem:3.2} \ A presentation of $\Omega^2(k)$ is determined completely
	by specifying four  elements $\rho_i$ of $\La$ 
	such that
\begin{align*} x\sigma + \rho_1\psi &=0, && y\sigma + \rho_2\psi =0\cr
	x\theta + \rho_3\psi&=0, && y\theta + \rho_4\psi =0.
\end{align*}
\end{lemma}

\bigskip

The module $\Omega^3(k)$ has  has dimension $2N-1$.\\
We take the generators $\sigma, \psi, \theta$ of $\Omega^2(k)$ in this order.
Then we identify
$$\Omega^3(k) = \{ (a_i) \in A^3\mid a_1\sigma  + a_2\psi  + a_3\theta  = 0\}.
$$
Using the  identities in Lemma \ref{lem:3.2}, we identify  the following four elements in $\Omega^3(k)$,
\begin{align*} f_1^{(3)} &= (x, \rho_1, 0)\cr
        f_2^{(3)} &=  (y, \rho_2, 0)\cr
        f_3^{(3)} &= (0, \rho_3, x)\cr
        f_4^{(3)} &= (0, \rho_4, y).
\end{align*}
The aim is to show that these are (independent) generators of $\Omega^3(k)$, by using only the conditions
(1) and (2).

\bigskip

\begin{lemma} \label{lem:3.3} We have $\Omega^3(k) = \sum_{i=1}^4 \La f_i^{(3)}$.
\end{lemma}

\bigskip

{\it Proof} \ In the proof we omit the
label $(-)^{(3)}$. 
Let $U_1= \La f_1 + \La f_2$ and let $U_2 = \La f_3 + \La f_4$, then $U_1+U_2\subseteq \Omega^3(k) \subseteq \La^3$. Let $\gamma_i$ be the projection
from $\Omega^3(k)$ onto the $i$-th component for $1\leq i\leq 3$.

(a) The image  $\gamma_1(U_1)$ is equal to  $\La x + \La y= \Omega(k)$, hence  has dimension $N -1$. 
Similarly  $\gamma_3(U_2) = \Omega(k)$ of dimension $N-1$.
Directly, the intersection $U_1\cap U_2$ is contained in $(0, \La, 0)$ and hence is in the kernel of $\gamma_1$ and $\gamma_3$. We have
\begin{align*}
\dim (U_1 + U_2)
&\leq \dim \gamma_1(U_1) + \dim \gamma_3(U_2)\cr
&\leq 2N-2
\end{align*}
which is one less than $\dim \Omega^3(k)$. 
Assume for a contradiction that  $U_1\cap U_2=0$, then $U_1+U_2=U_1\oplus U_2$, a  maximal submodule of $\Omega^3(k)$, which is indecomposable 
(since $\La$ is selfinjective).  As well ${\rm soc} \Omega^3(k)\cong {\rm top }\Omega^2(k)$ and is 3-dimensional. So the socle of $\La^3$ must be contained in $\Omega^3(k)$
 and hence there is an element 
$m:=(0, \omega, 0)\in \Omega^3(k)$ with $\omega$ an element spanning the socle of $\La$.  It follows that $\Omega^3(k) = U_1\oplus U_2 \oplus \langle m\rangle$, and is decomposable, a contradiction.
This shows that $U_1\cap U_2$ is non-zero and then must be 1-dimensional and 
spanned by $(0, \omega, 0)$ with $\omega$ in the socle of $\La$. 
The Lemma follows.
$\Box$

\bigskip

\begin{remark}\label{rem:3.3} \normalfont  With the notation as in the proof, the kernel of $\gamma_1$ is
	$$\{ (0, a_1\rho_1 + a_2\rho_2, 0)\mid (a_1, a_2) \in \Omega^2(k)\}$$
and similarly we can write down the kernel of $\gamma_2$. Hence we have
that if $\rho = {\rho_1\choose \rho_2}$  or ${\rho_3\choose \rho_4}$ then
$\sigma\cdot \rho$,  $\psi\cdot \rho$  \ and  $\theta\cdot \rho$ lie in the socle of $\La$. 
So far we  only use the information from (1) and (2).
In the two main types of algebras we want to deal with, we have more information. We make an assumption (which possibly is redundant).
\end{remark}

\begin{ass} \label{ass:3.4} We now give the third condition on the algebra $\La$, in addition to 
the conditions in Assumption \ref{ass:2.1}. 
We assume that there are elements $\rho_i\in \La$ satisfying the conditions in Lemma \ref{lem:3.2} and in addition 
\begin{align*} \sigma_1\rho_1+\sigma_2\rho_2&=0 && \psi_1\rho_3 + \psi_2\rho_4 = 0 \cr
        \psi_1\rho_1 + \psi_2\rho_2&=0 &&  \theta_1\rho_3 + \theta_2\rho_4 =0\cr
        \theta_1\rho_1+\theta_2\rho_2 &= \omega&&
	\end{align*}
where $0\neq \omega\in {\rm soc }\La$, and 
$$\theta_1\rho_1 + \theta_2\rho_2 + c(\sigma_1\rho_3+\sigma_2\rho_4)=0$$
for some $0\neq c\in k$. 
\end{ass}

\bigskip

\begin{expl}\label{ex:3.5} \normalfont Let $\La = A$ as in example \ref{ex:2.2}. 
We may take 
$\rho_1=\rho_4=0$ and $\rho_2 = q^{-n}x^{n-1}, \ \rho_3 =  -q^{m-1}y^{m-1}.$
Then  \label{ass:3.4} is satisfied:
\end{expl} 
\bigskip

\begin{lemma}\label{lem:3.5} Assume $\La = A$. Then:\\ 
	$$a_1\rho_1 + a_2\rho_2 = \left\{\begin{array}{ll} 0 & (a_1, a_2)=\sigma\cr
		                                      0&  (a_1, a_2)=\psi\cr
	q^{-n}y^{m-1}x^{n-1} & (a_1, a_2)=\theta. \end{array}\right.
$$
	$$a_3\rho_3 + a_4\rho_4 = \left\{\begin{array}{ll} -q^{m-1}x^{n-1}y^{m-1}  & (a_3, a_4)=\sigma\cr
                                                      0&  (a_3, a_3)=\psi\cr
        0 & (a_3, a_4)=\theta. \end{array}\right.
$$
\end{lemma}

We summarize the identities which have obtained. In the following
these are all we will use.

\begin{itemize}
	\item[(1)]
\begin{align*}
	\sigma_1x + \sigma_2y&=0\cr
	\psi_1x + \psi_2y &=0\cr
	\theta_1x + \theta_2y&=0
\end{align*}
\item[(2)]
\begin{align*} 
x\sigma_1 + \rho_1\psi_1&=0 && x\theta_1 + \rho_3\psi_1=0\cr
	x\sigma_2 + \rho_1\psi_2&=0 && x\theta_2 + \rho_3\psi_2=0\cr
&&&\cr
	y\sigma_1 + \rho_2\psi_1&=0 && y\theta_1+\rho_4\psi_1=0\cr
	y\sigma_2 + \rho_2\psi_2&=0 && y\theta_2 + \rho_4\psi_2=0
\end{align*}
\item[(3)]
	
	\begin{align*} \sigma_1\rho_1+\sigma_2\rho_2&=0 && \sigma_1\rho_3 + \sigma_2\rho_4 = (-c)\omega \cr
	\psi_1\rho_1 + \psi_2\rho_2&=0 &&  \psi_1\rho_3 + \psi_2\rho_4 =0\cr
	\theta_1\rho_1 + \theta_2\rho_2&=\omega && \theta_1\rho_3 + \theta_2\rho_4 =0
\end{align*}
\end{itemize}
Here $\omega$ is non-zero in the socle of $\La$, and $c$ is a non-zero scalar, and we have
$$\theta_1\rho_1 + \theta_2\rho_2 + c\sigma_1\rho_3 + c\sigma_2\rho_4=0.
$$

\bigskip

\section{A minimal projective resolution}\label{sec:4}

As before, we use left $\La$-modules, and we write maps to the right and compose from left to right. We will construct an exact sequence
$$\cdots \to P_{r+1}\stackrel{d_{r+1}}\to P_r\stackrel{d_r}\to \ldots \stackrel{d_3}\to P_2\stackrel{d_2}\to P_1\stackrel{d_1}\to P_0\to k \to 0
\leqno{(\dagger)}$$
with $P_r$ projective so that the image of $d_r$ is $\Omega^r(k)$. 

We define $P_r:= \La^{r+1}$
(for $r=0, 1, \ldots$). Then we define maps $d_r$ as matrices, where the matrix $d_r$ has $r+1$ rows and $r$ columns, 
and the action is given by applying $d_r$ to the standard generators of $\La^{r+1}$, and express the answer in terms of the standard 
generators of $\La ^r$.

First by the results of Section \ref{sec:3} we may take
$$d_1 = \left(\begin{matrix} x\cr y\end{matrix}\right), \ \  d_2 = \left(\begin{matrix} \sigma_1 & \sigma_2\cr \psi_1 & \psi_2\cr \theta_1 & \theta_2\end{matrix}\right),
d_3 = \left(\begin{matrix} x & \rho_1 &  \cr
y & \rho_2 & \cr
 & \rho_3&x\cr &\rho_4&y\end{matrix}\right), 
 \  d_4 = \left(\begin{matrix} \sigma_1 & \sigma_2 & & \cr \psi_1 & \psi_2 & & \cr
 \theta_1 & \theta_2 & c\sigma_1 & c\sigma_2\cr
 & & \psi_1 & \psi_2\cr
 & &\theta_1 & \theta_2\end{matrix}\right).
 $$
 so that the image of $d_1$ is precisely the radical of $\La$, as a submodule of $P_0$,  the image of $d_2$ is $\Omega^2(k)$, a submodule of $\La^2$, and so on.

\bigskip

\subsection{The matrix $d_r$ for $r\geq 5$ and $r$ odd}\ 
We define this recursively as follows. Let $r=2t-1$. 
Then take for $d_r$  the matrix with $r+1$ rows and $r$ columns, which we write in block form, 
$$d_{r} = \left(\begin{matrix}  d_{r-2}&  S_r\cr 0 & T\end{matrix}\right)
$$
Here $S_r$ is the matrix with $r-1$ rows and two columns, given as 
$$ S_r:= \left(\begin{matrix} S_r'\cr S_r''\end{matrix}\right), \ \ 
	S_r'=  \ 0_{(r-3)\times 2}, \ \ S_r''= \left(\begin{matrix}
			      \rho_1c^{t-2} & 0 \cr
\rho_2c^{t-2} & 0\end{matrix}\right).
$$
(That is,  $S_r'$ is the zero matrix of size $(r-3)\times 2$). 
Moreover, 
$T$ is the $2\times 2$ matrix 
$$T = \left(\begin{matrix} \rho_3 & x \cr \rho_4 & y\end{matrix}\right).
$$
\bigskip

\subsection{The matrix $d_r$ for $r\geq 4$ and $r$ even} \ 
Let $r=2t$. 
We define $d_r$ to be  the matrix with $r+1$ rows and $r$ columns which we write in block form:
$$d_r = \left(\begin{matrix} d_{r-2} &  U_r \cr
0 & V\end{matrix}\right).
$$
Here $U_r$ is the matrix with $r-1$ rows and two columns, given as 
$$U_r:= \left(\begin{matrix} U_r'\cr U_r''\end{matrix}\right), \ \  U_r'= 0_{(r-3)\times 2}, \ \ U_r''=
	\left(\begin{matrix} 0 & 0\cr 
	c^{t-1}\sigma_1 & c^{t-1}\sigma_2\end{matrix}\right).$$
Moreover,  $V$ is the matrix	
$$V = \left(\begin{matrix} \psi_1 & \psi_2 \cr
\theta_1 & \theta_2\end{matrix}\right).
$$

\bigskip

We show first that $d_{r+1}d_r = 0$ for each $r$. This will show that the above sequence
is a complex. Then we prove exactness.

\bigskip

\begin{lemma} \label{lem:4.1} The composition $d_{r+1}d_r$ is zero for each $r\geq 1$.
\end{lemma}

\bigskip

{\it Proof } 
For $r=1, 2, 3$ one checks this directly.
We continue by induction on $r$.

Assume first that $r$ is odd, and $2t-1 = r\geq 5$. Assume the
lemma holds for odd $s$ smaller than $r$. 
We compute $d_{r+1} d_r$ 
using the block shape, this is  equal to 
$$\left(\begin{matrix} d_{r-1}d_{r-2} & d_{r-1}S_r + U_{r+1}T\cr
0 & VT\end{matrix}\right).
$$ 
One checks that $VT=0$, using identities from (3) and (1). 
By the inductive hypothesis, $d_{r-1}d_{r-2} = 0$. It remains to show
that $d_{r-1}S_r + U_{r+1}T=0$. This is a matrix with $r$ rows and
two columns.
First,
$$d_{r-1}S_r = \left(\begin{matrix} d_{r-3} & U_{r-1}\cr 0 & V\end{matrix}\right)\cdot \left(\begin{matrix}S_r'\cr S_r''\end{matrix}\right) 
	= \left(\begin{matrix} U_{r-1}S_r''\cr 0 & VS_r''\end{matrix}\right).
$$
We find 
$$VS_r'' = \left(\begin{matrix} 0 & 0 \cr c^{t-2}(\theta_1\rho_1 + \theta_2\rho_2) & 0\end{matrix}\right).
$$
\bigskip

Next, 
$$U_{r-1}S_r'' = \left(\begin{matrix} U_{r-1}'\cr U_{r-1}''\end{matrix}\right)
	\cdot S_r'' = \left(\begin{matrix} 0 \cr
	U_{r-1}''S_r''\end{matrix}\right) = 0
$$
using that $\sigma_1\rho_1 + \sigma_2\rho_2=0$.
Next, 
$$U_{r+1}T = \left(\begin{matrix} U_{r+1}'\cr U_{r+1}''\end{matrix}\right)T
	= \left(\begin{matrix} 0 \cr
	U_{r+1}''T \end{matrix}\right)
	$$
and 
$$U_{r+1}''T = \left(\begin{matrix} 0 & 0 \cr c^{t-1}(\sigma_1\rho_3 + \sigma_2\rho_4) & 0\end{matrix}\right)
$$
We deduce that $d_{r-1}S_r + U_{r+1}T$ is the matrix of size $r\times 2$ 
with all entries zero, except possibly the $r1$, bu this is  equal to 
$$c^{t-2}(\theta_1\rho_1 + \theta_2\rho_2 + c\sigma_1\rho_3 + c\sigma_2\rho_4)=0
$$
Hence $d_{r+1}d_r=0$ for $r$ odd.

\bigskip

Now assume $r=2t \geq 6$, and assume true for even $s$ smaller than $r$. 
Then $d_{r+1} d_r$ is given by
$$\left(\begin{matrix} d_{r-1}d_{r-2} & d_{r-1}U_r + S_{r+1}V\cr 0 & TV\end{matrix}\right).
$$
One checks that $TV=0$, and by the inductive hypothesis also $d_{r-1}d_{r-2}=0$.
We show now that $d_{r-1}U_r + S_{r+1}V=0$. This is a matrix with $r$ rows and
two columns. We have
$$d_{r-1}U_r = \left(\begin{matrix} d_{r-3}& S_{r-1}\cr 0 & T\end{matrix}\right)\cdot \left(\begin{matrix}U_r'\cr U_r''\end{matrix}\right) 
	= \left(\begin{matrix}S_{r-1}U_r''\cr TU_r''\end{matrix}\right)
$$
We calculate
$$TU_r'' = \left(\begin{matrix} c^{t-1}x\sigma_1 & c^{t-1}x\sigma_2\cr
c^{t-1}y\sigma_1 & c^{t-1}y\sigma_2\end{matrix}\right).
$$
The matrix $S_{r-1}U_r''$ is zero. Next, we compute
$S_{r+1}V$, the first $r-2$ rows are zero and the last two rows are
$$\left(\begin{matrix} \rho_1\psi_1c^{t-1} & \rho_1 \psi_2c^{t-1}\cr
\rho_2\psi_1c^{t-1} & \rho_2\psi_2c^{t-1}\end{matrix}\right).
$$
Hence $d_{r-1}U_r + S_{r+1}V$ is the $r\times 2$ matrix with $r-2$ rows of
zeros, and where the last two rows are
$$c^{t-1}\left(\begin{matrix} x\sigma_1 + \rho_1\psi_1 & x\sigma_2 + \rho_1\psi_1\cr
y\sigma_1 + \rho_2\psi_1 & y\sigma_2 + \rho_2\psi_2\end{matrix}\right).
$$
This is zero, by the identities at the end of Section \ref{sec:3}.
$\Box$

\bigskip

This shows that $(\dagger)$ is a complex, and next we will show that
it is exact.

\bigskip

\begin{prop}Let $r= 2t$ or $2t-1$ for $t\geq 1$. 
	We have ${\rm im} ( d_r) = \Omega^r(k)$ of dimension
	 $tN+1$ if $r$ is even, or $tN-1$ if $r$ is odd. 
	
\end{prop}

\bigskip

{\it Proof} This is true for $r=1, 2, 3$ by Section \ref{sec:3}, and let $r =2t$ or $2t-1$ where $r\geq 4$. As an inductive
hypothesis, assume true for all $m\leq r-1$.  Then 
${\rm im} (d_{r-1}) = \Omega^{r-1}(k)$ of dimension $tN-1$ or $(t-1)N+1$, depending
on the parity of $r$.  Since $P_{r-1} = \La^r$, it  follows that $\Omega^r(k)$ has
dimension $tN+1$ or $tN - 1$; and moreover we know
${\rm Im} (d_r)\subseteq \Omega^r(k)$ since $d_rd_{r-1}=0$.
We are done if we show that ${\rm Im} (d_r)$ has 
the same dimension as $\Omega^r(k)$. 

\bigskip
(1) Assume first $r$ is even, so that  $r=2t \geq 4.$
We must show that ${\rm im} (d_r)$ has dimension $tN+1$. By the above, the dimension
is at most $tN+1$. 
 Let $\pi: \La^{r+1}\to \La^{r-1}$ be the projection onto the first
 $r-1$ coordinates. This takes $f_i^{r}$ to $f_i^{r-2}$ for
 $1\leq i\leq r-1$, and it takes $f_r^r$ and $f_{r+1}^r$ to zero.
 Hence $ ({\rm Im}(d_r))\pi = {\rm Im } (d_{r-2})$ which has
 dimension $(t-1)N+1$. Furthermore, 
 the submodule $\langle f_r^{r}, f_{r+1}^r\rangle$ contained in the kernel
 of $\pi$ is isomorphic to $\langle \psi, \theta\rangle$ and has
 dimension $N$ (see the remark below \ref{lem:3.1}). So the dimension of ${\rm Im} (d_r)$ is 
 at least $N + (t-1)N + 1$, and then equality follows.
 
 \bigskip

 (2) Now assume $r$ is odd so that $r=2t-1$, we must show that
 ${\rm im} (d_r)$ has dimension $tN-1$. By the above, the dimension is at most
 $tN-1$.
Write ${\rm Im}(d_r) = U + V\subseteq \La^r$ where
$U$ is the submodule generated by $f_1^r, \cdots, f_{r-1}^r$ and 
$V$ is the submodule generated by $f_r^r, f_{r+1}^r$. 

We find the dimension of $V$.  Let 
$\pi: \La^r\to \La$ be the projection onto the last coordinate. 
Then $(V)\pi = \La x + \La y$ and hence has dimension $N-1$.
We consider the kernel of $\pi$ restricted to $V$.  As a vector space it is 
isomorphic to
$$\{ a_1\rho_3 + a_2\rho_4\mid (a_1, a_2) \in \Omega^2(k)\}.$$
This is computed in Section \ref{sec:3}, and it is just the socle of $\La$
(concentrated in  the $r-1$-th coordinate), it
is 1-dimensional, and hence $\dim V = N$. 

Now we find the dimension of $U$. Let $\pi': \La^r\to \La^{r-2}$ be the projection onto the first $r-2$ components. 
We have $(f_i^r)\pi' = f_i^{r-2}$ for $1\leq i\leq r-1$. 
So $(U)\pi' = {\rm Im} (d_{r-2})$ of dimension $(t-1)N-1$, and
hence $\dim U$ is at least $(t-1)N-1$. 
We can identify the kernel of $\pi'$ restricted to $U$. 
Note that $\pi'$ restricted to $\langle f_i^r\mid 1\leq i\leq r-4 \rangle$
is injective, we find that 
${\rm ker}(\pi'|_U)$ is isomorphic to
$$\{ a_1\rho_1 + a_2\rho_2\mid (a_1, a_2) \in \Omega^2(k)\}.$$
This  is again the socle of $\La$ (in the $r-1$-th coordinate)
and is 1-dimensional. So $\dim U = (t-1)N$. Moreover, we see
that $U\cap V$ is 1-dimensional and then $\dim (U+V) = tN-1$ as required.
$\Box$

\vspace*{1cm}

\section{The cohomology ring}\label{sec:5}

By Section \ref{sec:3}, for each $r\geq 1$, the module $\Omega^r(k)\subseteq \La^{r}$ 
is generated by the $r+1$ rows of the matrix $d_r$. Write $f_i^r$ for the
$i$-th row of this matrix, for $1\leq i\leq r+1$. 
Then the $f_i^r$ form a minimal set of generators, and hence we can take
as a vector space basis for $H^r(\La, k) = {\rm Hom}_{\La}\Omega^r(k), k)$ 
the set $\{ \vf_i^r\}$ which is the dual basis for the $f_i^r$.

If $f, g$ are homogeneous elements of degrees
$s, t$ (say), then we take the product $f\cdot g$ as the class of
$\Omega^t(f)\cdot g$ (recall that we apply maps to the right).

We will show  that for $t\geq 1$ we have
\begin{align*} H^{2t}(\La, k)\times H^1(\La, k) &= H^{2t+1}(\La, k), \cr
H^{2t}(\La, k)\times H^2(\La, k) &= H^{2t+2}(\La, k).
\end{align*}

In this section we write just $H^m$ meaning $H^m(\La, k)$.

\bigskip

\begin{lemma}\label{lem:5.1}  We have $H^{2t}\times H^1 = H^{2t+1}$.
More precise, for $0\leq s\leq t$, the following hold
$$\vf_{2s+1}^{2t}\cdot \vf_1^1 = \vf_{2s+1}^{2t+1}, \ \ 
	\vf_{2s+1}^{2t}\cdot \vf_2^1 = \vf_{2s+2}^{2t+1}.
$$
That is, $H^{2t}\times H^1$ contains a basis of $H^{2t+1}$

\end{lemma}

\bigskip

{\it Proof } 
The space $H^{2t}(\La, k)$ has basis
$\{ \vf_i^{2t} \mid 1\leq i\leq 2t+1\}$. 
The product $\vf_i^{2t}\cdot \vf_j^1$ is then the class of
$$\Omega(\vf_i^{2t}) \vf_j^1 \ \ (1\leq j\leq 2).
$$
The diagram to compute $\Omega(\vf_i^{2t})'s$ is of the form
$$\CD  \Omega^{2t+1}(k) @>>> P_{2t} @>{d_{2t}}>> \Omega^{2t}(k) @>>> 0\cr
@V{\Omega(\vf_i^{2t})}VV  @V{h_i}VV  @V{\vf_i^{2t}}VV \cr
\Omega(k) @>>> \La @>{d_0}>> k @>>> 0
\endCD
$$
The composition $d_{2t} \vf_i^{2t}$ is equal to $\ve_i^*$, the dual
of the standard basis element. We can therefore take
$$h_i = \ve_i^{T}
$$
(the $i$th element in the  standard basis of  column vectors). 
Then $\Omega(\vf_i^{2t})$ is the restriction of $h_i$ to $\Omega^{2t+1}(k)$, which
is the submodule of $\La^{2t+1}$ generated by the rows of $d_{2t+1}$, which 
we call $f_1^{2t+1},  f_2^{2t+1}, \cdots, f_{2t+2}^{2t+1}$. 

\bigskip
Let $i=2s+1$ for $0\leq s\leq t$. Then $h_i$ takes
\begin{align*} 
	f_i^{2t+1} &\mapsto  x\cr
	f_{i+1}^{2t+1} &\mapsto y\cr
	f_u^{2t+1} &\mapsto 0  (u\neq i, i+1)
\end{align*}
From this the claim follows.
$\Box$

\bigskip

Recall that the space $H^{2t}(\La, k)$ has basis 
$\{ \vf_i^{2t} \mid 1\leq i\leq 2t+1\}$.

\bigskip
\begin{lemma} \label{lem:5.2} We have $H^{2t}\times H^2 = H^{2t+2}$. 
More precisely, for $1\leq s\leq t$ we have
$$\vf_{2s+1}^{2t}\cdot \vf_1^2 =  c^s\vf_{2s+1}^{2t+2},  \ \vf_{2s+1}^{2t}\cdot \vf_2^2 = \vf_{2s+2}^{2t+2} \mbox{ \ and } \ \vf_{2s+1}^{2t}\cdot \vf_3^2 = \vf_{2s+3}^{2t+2}.
$$
Hence $H^{2t}\times H^2$ contains a basis for $H^{2t+2}$. 
\end{lemma}

\bigskip

{\it Proof}  Let $1\leq i\leq 2t+1$. 
The product $\vf_i^{2t}\cdot \vf_j^2$ is  the class of
$$\Omega^2(\vf_i^{2t})\circ \vf_j^2 \ \ (1\leq j\leq 3).
$$
We take as the diagram to compute $\Omega^2(\vf_i^{2t})'s$ the extension of 
the previous diagram, that is we take
$$\CD  \Omega^{2t+2}(k)  @>>> P_{2t+1} @>{d_{2t+1}}>> P_{2t} @>{d_{2t}}>> \Omega^{2t}(k) @>>> 0 \cr
@V{\Omega^2(\vf_i^{2t})}VV  @V{\ell_i}VV   @V{h_i}VV  @V{\vf_i^{2t}}VV \cr
 \Omega^2(k) @>>> \La^2 @>{d_1}>> \La @>{d_0}>> k @>>> 0
\endCD
$$
As before, we take $h_i:= \ve_i^T$ (the standard column vector). 
Now assume $i$ is an odd number $\leq 2t+1$,  write 
$i =2s+1 \ (0\leq s\leq t).$
Let $\ell_i$ be the $(2t+2)\times 2$ matrix 
whose row $i$ is $(1 \ 0)$ and whose row $i+1$ is $(0 \ 1)$, and where all other entries are zero. 
Then $\ell_i\cdot d_1$ is the column vector $x\ve_i + y\ve_{i+1}$ and this is the same as 
$d_{2t+1}\circ h_i$. So the map given by $\ell_i$ lifts $h_i$. Therefore
$\Omega^2(\vf_i^{2t})$ can be taken as the restriction of $\ell_i$ to $\Omega^{2t+2}(k)$, ie the submodule
of $\La^{2t+2}$ generated by the rows of the matrix $d_{2t+2}$.

Thus $\ell_i$ takes
\begin{align*} f_i^{2t+2} &\mapsto c^s\sigma, \cr
f_{i+1}^{2t+2} & \mapsto  \psi\cr
f_{i+3}^{2t+2} &\mapsto \theta\cr
f_{j}^{2t+2} &\mapsto 0  \ \mbox{else} 
\end{align*}
Recalling that $f_1^2= \sigma, f_2^2=\psi$ and $f_3^2= \theta$, we obtain
the statement of the Lemma.
$\Box$

\bigskip

\begin{theorem}\label{thm:5.3} The algebra $H^*(\La, k)$ is finitely generated, and is generated
by elements in degrees $1, 2$. 
\end{theorem}

\bigskip

This follows directly from the above two lemma, by induction on $t$. 
We note that this only uses little information on 
the algebra $\La$.  To get a presentation of the cohomology, 
requires further detailed information. In the following we give some illustrations.

\subsection{Products of two elements of degree 1}

We have $H^1(\La, k) = {\rm Hom}_{\La}(\Omega(k), k)$.  We have $f_1^1=x$ and $f_2^1=y$ as minimal generators
of $\Omega(k)$, and then  $\vf_i^1$ are corresponding dual basis elements of $H^1(\La, k)$.
We must compute $\Omega(\vf_i^1)$ for $i=1, 2$. 
Hence we need to find $h_i$ lifting $\vf_i^1$ and then $\Omega(\vf_i^1)$ is the restriction of $h_i$ to $\Omega^2(k)$.

The relevant diagram is
$$\CD  \Omega^2(k) @>>> P_1 =\La^2 @>{d_1}>> \Omega(k) @>>> 0\cr
@V{\Omega(\vf_i^1)}VV  @V{h_i}VV  @V{\vf_i^1}VV \cr
\Omega(k) @>>> \La @>{d_0}>> k @>>> 0
\endCD
$$
The composition $d_1\circ \vf_i^1$ is the dual $\ve_i^*$ of the standard  basis element  for $\La^2$. We can take
$h_i:= \ve_i^T$, the column vector with $1$ in place $i$ and zero elsewhere. 
Then $\Omega(\vf_i^1)$, the restriction of $h_i$ to $\Omega^2(k)$, takes each generator of
$\Omega^2(k)$ to its $i$th  coordinate. 
The  generators of $\Omega^2(k)$ are
$$f_1^2=\sigma, \ \ f_2^2 = \psi, \ \ f_3^2 = \theta.
$$
The $\sigma_i, \psi_i, \theta_i$ are in the radical of $\La$, so there are expressions
\begin{align*} \sigma_i =& u_{ix}x + u_{iy}y \cr
\psi_i =& v_{ix}x + v_{iy}y\cr
\theta_i =& w_{ix}x + w_{iy}y
\end{align*}
with coefficients in $\La$.  

Let $\ol{(-)}:  \La\to k = \La/{\rm rad}\La$ be the canonical map. Then 
for any $p\in {\rm rad} \La$ if we write $p = cx + dy$, then 
$$ \vf_1^1(p) = \ol{c} , \ \ \vf_2^1(p) = \ol{d}.
$$
Therefore 
\begin{align*}
\vf_1^1\cdot \vf_1^1 &=  \ol{u_{1x}} \vf_1^2 +  \ol{v_{1x}}\vf_2^2 + \ol{w_{1x}}\vf_3^2\cr
\vf_1^1\cdot \vf_2^1 &=  \ol{u_{1y}} \vf_1^2 +  \ol{v_{1y}}\vf_2^2 + \ol{w_{1y}}\vf_3^2\cr
\vf_2^1\cdot \vf_1^1 &=  \ol{u_{2x}} \vf_1^2 +  \ol{v_{2x}}\vf_2^2 + \ol{w_{2x}}\vf_3^2\cr
\vf_2^1\cdot \vf_2^1 &=  \ol{u_{2y}} \vf_1^2 +  \ol{v_{2y}}\vf_2^2 + \ol{w_{2y}}\vf_3^2.
\end{align*}

\bigskip

Since the $\psi_1, \psi_2$ are independent generators of the radical, and so are $x, y$, we would know that
the matrix $\left(\begin{matrix} \ol{v_{1x}} & \ol{v_{1y}}\cr \ol{v_{2x}} & \ol{v_{2y}}\end{matrix}\right)$ is invertible.
With this we would get $\vf_2^2$ in $H^2$, but we cannot say much more without knowing the algebra.

\begin{expl} \label{ex:5.1}  Take $\La = A$. Then 
$$\vf_1^1\cdot \vf_1^1 = \left\{\begin{array}{ll} \vf_1^2 & n = 2\cr 0 & \mbox{else}
\end{array}
\right. 
$$
Similarly $\vf_2^2\cdot \vf_2^2 = \vf_2^2$ if $m=2$, and is zero otherwise.
Furthermore, $\vf_1^2\cdot \vf_2^1= -q\vf_2^2$ and 
$\vf_2^1\cdot \vf_1^1 = \vf_2^2$.
\end{expl}

\bigskip

\parindent0pt

\subsection{The products  $H^2\times H^2$.}

In Lemma \ref{lem:5.2}  we  have dealt with the products  $\vf_i^2\cdot \vf_j^2$ where $i$ is odd. 
We consider   products with $i=2$, so we must find a suitable map $\ell_2$ such that
$\ell_2 \cdot d_1 = d_3\cdot h_2$.

\medskip
The map $d_3\cdot h_2$ takes $\ve_j\mapsto \rho_j$ for $1\leq j\leq 4$.
Since $\rho_j$ is in the radical of $\La$, there are elements $u_j,  v_j\in \La$ such that
$\rho_j = u_jx + v_jy.$  We fix such elements, and define
$$\ell_2:= \left(\begin{matrix} u_1&v_1\cr u_2 & v_2\cr u_3 & v_3\cr u_4 & v_4\end{matrix}\right).
$$
Then $\ell_2\cdot d_1=  d_3\cdot h_2$. 
Now
$\Omega^2(\vf_2^2)$ is the restriction of $\ell_2$ to the submodule 
$\langle f_j^4\mid 1\leq i\leq 5\rangle$  of $\La^4$. 
The elements $f_i^4\cdot \ell_2$  must be  $\La$-combinations of $\sigma, \psi, \theta$. 
In general, no easy description can be expected.
We will discuss this for the algebra $\La = A5$ in the next section.

\medskip

We observe that in general the cohomology is not graded commutative, in
fact if $c\neq 1$ then the even part is not commutative.

\newpage

\section{The Hopf algebra $A5$ for $p\geq 3$ with $\beta \neq 0$}\label{sec:6}

We consider the algebra, $A5(\beta)$, as in the Introduction. Here $z$ plays the role of $x$ (and $y$ is kept). We have
$\La = k\langle z,y \rangle/I$, and we have the following.
We set $a:=yz - zy$, it is central, and with this, the ideal $I$ is determined by
$$y^p=0, \ a^p=0, \ z^p + a^{p-1}y - \beta a = 0.
$$
We use freely that $a$ is central and $a^p=0$.
In particular we will use the fact that $a^{p-1}yz = a^{p-1}zy.$

\medskip

We define 
$$\sigma:= (z^{p^2-1}, 0)  \ \ \psi = (z^{p-1} - \beta y, a^{p-1} + \beta z), \ \ \theta = (0, y^{p-1}).
$$
We will show that these satisfy Assumptions \ref{ass:2.1} and \ref{ass:3.4}.

\medskip

\noindent We start with some identities. By definition $z^p= \beta a - a^{p-1}y$. 
The first observation is that
$$z^pa^{p-1}=0.$$
Next, $(z^p)^2 = \beta^2a^2$ and hence $(z^p)^r= \beta^ra^r$ for $2\leq r$ and 
$z^{p^2}=0$. 
Moreover
$$z^{p^2-2} = \beta^{p-1}a^{p-1}z^{p-2} \ \mbox{ and } \ z^{p^2-1} = (z^p)^{p-1}z^{p-1} = \beta^{p-1}a^{p-1}z^{p-1} \neq 0.
$$
Hence $z^{p^2-1}a=0$ and $z^{p^2-1}$ commutes with $y$.
For computations, we  will use the following
formula which allows one to interchange powers of $y$ with powers of $z$. This  can
be proved by induction from $a = yz-zy$.
\bigskip

\begin{lemma}\label{lem:6.1}
We have
$$y^bz^c = \sum_{r\leq b, c} r! {b\choose r}{c\choose r}a^rz^{c-r}y^{b-r}$$
\end{lemma}

\bigskip

Next, we describe a basis for the algebra.

\begin{lemma} \label{lem:6.2} The algebra $\La$
has $k$-basis
$$\{ a^iy^jz^k \mid 0\leq i, j, k\leq p-1\}.$$
The element $a^iy^jz^k$ in this set belongs to $J^m\setminus J^{m+1}$ where
        $m= pi + j + k$. The socle of $A$ is spanned
        by $a^{p-1}y^{p-1}z^{p-1}$.
\end{lemma}

This is straightforward. Next, we show that the Assumption \ref{ass:2.1} is satisfied.

\bigskip

\begin{lemma} \label{lem:6.3}
The elements $\sigma$ and $\theta$ are not in $A\psi$, and are linearly independent.
\end{lemma}

\bigskip

{\it Proof } Clearly $\sigma$ and $\theta$ are linearly independent over $k$.  Assume for a contradiction that
$\sigma$ belongs to $A\psi$. That is, there is some $\gamma \in A$ such that
$\sigma = \gamma \psi$.
That is
$$z^{p^2-1} = \gamma(z^{p-1}-\beta y), \ \ 0 = \gamma(a^{p-1}+\beta z).
$$
Consider the first identity,
$z^{p-1}-\beta y \in J\setminus J^1$ and the product is in
$J^{p^2-1}$. This means that $\gamma$ must be in $J^{p^2-2}\setminus J^{p^2-3}$. 
A basis for this quotient consists
of the cosets of elements
$$a^iy^jz^k, \ \ pi + j+k = p^2-2.$$
Solving this, shows that $i=1$ or $2$ and
$\gamma$ has highest order term
$$a^{p-1}y^jz^k, \ \ j+k=p-2 \ \ \mbox{ or } \ a^{p-2}y^{p-1}z^{p-1}.
$$
In both cases there is no solution for
$\gamma(z^{p-1}-\beta y) = \beta^{p-1}(a^{p-1}z^{p-1})$, a contradiction.

Next, assume for a contradiction that there is $\delta\in A$ such that $\theta = \delta \psi$.
That is
$$0 = \delta\psi_1 = \delta(z^{p-1}-\beta y), \ \  y^{p-1} = \delta(a^{p-1}+\beta z).
$$
Consider the second identity.  Modulo the ideal $I:= a^{p-1}A$
we require  $\beta\delta z = y^{p-1}$, which has
no solution, a contradiction.
$\Box$

\bigskip

\begin{lemma} \label{lem:6.4} Let
\begin{align*} 
\rho_1&=0= \rho_4,\cr
\rho_2 &= \beta^{-1}z^{p^2-1}\cr
\rho_3 &= -\beta^{-1}y^{p-1} - \beta^{-2}(y^{p-2}z^{p-1}) + \beta^{-3}(a^{p-1}y^{p-2}z^{p-2}).
\end{align*}
These satisfy Assumption \ref{ass:3.4}.
\end{lemma}

\bigskip

{\it Proof} \ 
It is clear that $\rho_1$ and $\rho_4$ satisfy Lemma \ref{lem:3.2} \\
(1) To show that $\rho_2$ satisfies Lemma \ref{lem:3.2}, we need

(i) $yz^{p^2-1} + \beta^{-1}z^{p^2-1}(z^{p-1}-\beta y) = 0$ and\\
(ii) $\beta^{-1}z^{p^2-1}(a^{p-1} + \beta z) = 0$.

For (i), we have $z^{p^2-1}z^{p-1}=0$ and hence (i) becomes 
$$yz^{p^2-1} - z^{p^2-1}y = 0.
$$
To prove (ii), we have $z^{p^2-1}a=0$ and hence (ii) becomes $z^{p^2}=0$.

\medskip

(2) We show now that $\rho_3$ satisfies Lemma \ref{lem:3.2}, that is we must show\\
(i) $\rho_3\psi_1=0$ and\\
(ii) $zy^{p-1} + \rho_3\psi_2=0$.

(i) \ We compute the product of each of the
three terms of $\rho_3$  with $\psi_1$, ie with $z^{p-1} - \beta y$.

$$-\beta^{-1}y^{p-1}(z^{p-1}-\beta y) = -\beta^{-1}y^{p-1}z^{p-1}
$$
To compute the second product, observe
$$z^{p-1}z^{p-1} = (\beta a - a^{p-1}y)z^{p-2} \ \mbox{ and } \ 
z^{p-1}y = yz^{p-1} + az^{p-2}
$$
So 
\begin{align*} -\beta^{-1}y^{p-2}z^{p-1}(z^{p-1}-\beta y)=&
	-\beta^{-2}y^{p-2}(\beta a-a^{p-1}y)z^{p-2}
	+ \beta^{-1}y^{p-2}(yz^{p-1} + az^{p-2})\cr
	=& -\beta^{-1}ay^{p-2}z^{p-2} + \beta^{-2}a^{p-1}y^{p-1}z^{p-2} 
	+ \beta^{-1}y^{p-1}z^{p-1} + \beta^{-1}ay^{p-2}z^{p-2}\cr
	=& \beta^{-2}a^{p-1}y^{p-1}z^{p-2} + \beta^{-1}y^{p-1}z^{p-1}.
\end{align*}

For the third product,
$$\beta^{-3}a^{p-1}y^{p-2}z^{p-2}(z^{p-1}-\beta y) = 0 - \beta^{-2}a^{p-1}y^{p-1}z^{p-2}.
$$
The sum of the three terms we obtained is equal to zero, as stated.

(ii) We compute each of the three terms of $\rho_3\psi_2$.
First
$$-\beta^{-1}y^{p-1}\psi_2= -\beta^{-1}a^{p-1}y^{p-1} - y^{p-1}z$$
For the next product, we use $y^{p-2}z^p = \beta ay^{p-2} - a^{p-1}y^{p-1}$ and get
\begin{align*}
	-\beta^{-2}y^{p-2}z^{p-1}\psi_2=&
	-\beta^{-2}a^{p-1}y^{p-2}z^{p-1} - ay^{p-2} + \beta^{-1}a^{p-1}y^{p-1}
\end{align*}
Finally 
\begin{align*} \beta^{-3}a^{p-1}y^{p-2}z^{p-2}\psi_2&=
	\beta^{-3}a^{p-1}y^{p-2}z^{p-2}a^{p-2} + \beta^{-2}a^{p-1}y^{p-2}z^{p-1}\cr
	&= \beta^{-2}a^{p-1}y^{p-2}z^{p-1}.
\end{align*}
Adding these we see
$$\rho_3\psi_2 = -y^{p-1}z - ay^{p-2}.
$$
Using the formula \ref{lem:6.1} we have
$$y^{p-1}z = zy^{p-1} - ay^{p-2}
$$
and $\rho_3\psi_2 + zy^{p-1}=0$ as required.
$\Box$

\bigskip

We have proved that Assumption \ref{ass:2.1} holds. It remains to verify the identities
in Assumption \ref{ass:3.4}.

Three of these are easy,  
\begin{align*} \sigma_1\rho_1 + \sigma_2\rho_2&=0\cr
\theta_1\rho_3+\theta_2\rho_4&=0\cr
\psi_1\rho_1 + \psi_2\rho_2&= (a^{p-1}+\beta z)\beta^{-1}z^{p^2-1} = 0.
\end{align*}
	As well
$$\theta_1\rho_1+\theta_2\rho_2 = \beta^{-1}y^{p-1}z^{p^2-1} = 
	\beta^{p-2}(a^{p-1}y^{p-1}z^{p-1})
$$
which is a non-zero element in the socle of $\La$. Moreover
$$
\sigma_1\rho_3 + \sigma_4\rho_4 =  
	\beta^{p-1}a^{p-1}z^{p-1}\rho_3 = -\beta^{p-2}a^{p-1}y^{p-1}z^{p-1}
$$
where the last equality follows by noting that
$(a^{p-1}z^{p-1})$ is annihilated by $a$ and by $z$
and that $z, y$ commute on elements which have a factor $a^{p-1}$.
	This gives the last equation in Assumption \ref{ass:3.4}, where the scalar $c$ is equal to
	$1$. 
	
\bigskip

It remains to show that $\psi_1\rho_3=0$.
We compute each of the three terms. First
$$\psi_1\rho_3 = (z^{p-1}-\beta y)(-\beta^{-1}y^{p-1} )
	= -\beta^{-1}(z^{p-1}y^{p-1}) \leqno{(\dagger)}
$$
The second product is equal to 
	$$\beta^{-2}(z^{p-1}y^{p-2}z^{p-1}) + \beta^{-1}(y^{p-1}z^{p-1}).
	$$
We write (I) for the first summand and (II) for the second. 
	Consider (I), this needs
\begin{align*}
	z^{p-1}y^{p-2}z^{p-1} =& z^{p-1}\sum_{r\leq p-2}r!{p-2\choose r}{p-1\choose r}a^rz^{p-1-r}y^{p-2-r}\cr
	=& z^p\sum_{r\leq p-2} r!{p-2\choose r}{p-1\choose r}a^rz^{p-2-r}y^{p-2-r}\cr
	=& \beta[\sum_{r\leq p-2}r!{p-2\choose r}{p-1\choose r}a^{r+1}z^{p-2-r}y^{p-2-r}] \ 
	- a^{p-1}yz^{p-2}y^{p-2} \cr
	=& \beta(\Sigma_1) \  - a^{p-1}z^{p-2}y^{p-1}
\end{align*}
defining $\Sigma_1$ to be the sum. 

	Similarly (II) is equal to
$$\beta^{-1}y^{p-1}z^{p-1} = \beta^{-1}(\sum_{r\leq p-1}r!{p-1\choose r}^2a^rz^{p-1-r}y^{p-1-r}) =: \beta^{-1}\Sigma_2.
$$
The third product is equal to
$$(z^{p-1}-\beta y)\beta^{-3}a^{p-1}y^{p-2}z^{p-2} = -\beta^{-2}a^{p-1}y^{p-1}z^{p-2}.
$$
In total we have
$$\psi_1\rho_3 = (\dagger) - \beta^{-2}\beta\Sigma_1 + \beta^{-2}a^{p-1}z^{p-2}y^{p-1}
+ \beta^{-1}\Sigma_2 - \beta^{-1}a^{p-1}y^{p-1}z^{p-2}.
$$
The term $(\dagger)$ cancels agains the term for $r=0$ in $\Sigma_2$. The two 
monomials with $\beta^{-2}$ cancel. This leaves 
$$\Delta:= -\beta^{-1}\Sigma_1 + \beta^{-1}\sum_{r=1^{p-1}}r!{p-1\choose r}^2a^rz^{p-1-r}y^{p-1-r}.
$$
Changing variables in the second sum, we can combine the terms and get
$$\Delta = \beta^{-1}\sum_{r=0}^{p-2} c_ra^{r+1}z^{p-2-r}y^{p-2-r}
$$
where
$$c_r = -r!{p-2\choose r}{p-1\choose r} + (r+1)!{p-1\choose r+1}^2.
$$
By an elementary calculation one gets that this is zero in $k$. 
This completes the proof that the algebra $A5$ satisfies the identities in \ref{ass:3.4}.
$\Box$

\bigskip

\subsection{Products in $H^2\times H^2$ for the algebra $A5$}

We will continue with the notation of  Lemma \ref{lem:5.2},  in particular it will confirm that
the even part of $H^*(\La, k)$ is indeed commutative (which is the case for cohomology rings of
Hopf algebras).

\bigskip

\begin{lemma} For each $1\leq i, j\leq 3$ 
	we have $\vf_i^2\cdot \vf_j^2 = \vf_j^2\cdot \vf_i^2$. 
\end{lemma}

{\it Proof } Most of this is proved in Lemma \ref{lem:5.2}. Note that we
do not need to deal with squares.
We are left to show that 
$$\vf_2^2\cdot \vf_1^2= \vf_2^4 \ \ \mbox{ and } 
\vf_2^2\cdot \vf_3^2 = \vf_4^4$$

Recall $\rho_1=\rho_4=0$ and we have
$$\ell_2 = \left(\begin{matrix} 0 & 0\cr u_2 & v_2\cr u_3 & v_3\cr 0 & 0\end{matrix}\right).$$
where $\rho_j = u_jz + v_jy$ for $j=2, 3$. 
Then 
for each $j$ we know that
$$f_j^4\cdot \ell_2 = \alpha_j\sigma + \beta_j\psi + \gamma_j\theta
$$
where $\alpha_j, \gamma_j \in k$ are unique 
(from Lemma \ref{lem:6.3}), and where $\beta_j \in A$. We will only need
to identify the $\alpha_j$ and $\gamma_j$.

First we have $f_1^4\cdot \ell_2=0$ and $f_5^4\cdot \ell_2 =0$. Moreover
we see directly that $f_2^4\cdot \ell_2 = \sigma$, so $\alpha_2=1$ and $\gamma_2=0$, as required. 

Next, we need 
$$f_3^4\cdot \ell_2 = (y^{p-1}(\beta^{-1}z^{p^2-2}) + z^{p^2-1}u_3, 0 + z^{p^2-1}v_3)
$$
Since $u_3$ and $v_3$ are in the radical, it is clear that 
there is no non-zero summand $\alpha_3\sigma$ or $\gamma_3\theta$. 
That is $\alpha_3=0$ and $\gamma_3=0$. 
Similarly we see $\gamma_4=1$ and $\alpha_4=0$.
$\Box$

\end{document}